\documentclass[12pt, a4paper]{article}
\newtheorem{thm}{{\bf T}{\footnotesize \bf HEOREM}}
\newtheorem{lm}[thm]{{\bf L}{\footnotesize \bf EMMA}}
\newtheorem{cor}[thm]{{\bf C}{\footnotesize \bf OROLLARY}}
\newtheorem{pro}[thm]{{\bf P}{\footnotesize \bf ROPOSITION}}
\newtheorem{conj}[thm]{{\bf C}{\footnotesize \bf ONJECTURE}}

\newcommand{\qed}{\hbox{\rule{6pt}{6pt}}}
\newcommand{\dg}{\mbox{deg}\,}
\newcommand{\dgg}{\mbox{deg}}
\newenvironment{Proof}{\noindent {\it Proof.} }{\hbox{\rule{6pt}{6pt}}
\bigskip}
\newenvironment{Proofof}[1]
{\noindent {\it Proof of {#1}.} }{\hbox{\rule{6pt}{6pt}} \bigskip}
\newcommand{\qqed}{\hbox{\rule{4pt}{5pt}}}
\usepackage[utf8]{inputenc}
\usepackage[dvipdfmx]{graphicx}
\usepackage{tikz}
\usepackage{tikzinclude}
\usepackage{amsmath, amsfonts}
\usepackage{makecell}
\renewcommand{\baselinestretch}{1}

\begin{document}
% \linenumbers
\sloppy
\title{The matching extendability of optimal $1$-embedded graphs on the projective plane}
\author{
Shohei Koizumi\thanks{%
Graduate School of Science and Technology,
Niigata University,
8050 Ikarashi 2-no-cho, Nishi-ku, Niigata, 950-2181, Japan.
E-mail: {\tt s-koizumi@m.sc.niigata-u.ac.jp}} \
Yusuke Suzuki\thanks{%
Department of Mathematics, Niigata University, 
8050 Ikarashi 2-no-cho, Nishi-ku, Niigata, 950-2181, Japan.
Email: {\tt y-suzuki@math.sc.niigata-u.ac.jp}}
%\thanks{%
%Partly supported by Japan Society for the Promotion of Science,
%Grant-in-Aid for Scientific Research (C) 16K05250.}
}
\date{}
\maketitle

\begin{abstract}
\noindent
In this paper, we discuss matching extendability 
of optimal $1$-projective plane graphs (abbreviated as O1PPG), 
which are drawn on the projective plane $P^2$ so that 
every edge crosses another edge at most once, and has $n$ vertices and exactly $4n- 4$ edges.
We first show that every O1PPG of even order is $1$-extendable. 
Next, we characterize $2$-extendable O1PPG's in terms of 
a separating cycle consisting of only non-crossing edges. 
Moreover, we characterize O1PPG's having connectivity exactly $5$.  
Using the characterization, we further identify three independent edges in those graphs 
that are not extendable. 
\end{abstract}

\noindent
{\bf Keywords:} 1-embeddable graph,  projective plane, perfect matching 

\noindent
{\bf 2020 Mathematical Subject Classification:} 05C10

\section{Introduction}\label{sect:intro}
\noindent
Our graphs dealt in this paper are all finite, simple and connected. 
%Furthermore, ``graph'' in this paper indicates ``simple graph'' 
%unless otherwise specified; only in one section, multigraphs make a brief appearance 
%in our proof. 
We denote the vertex set and the edge set of a graph $G$ 
by $V(G)$ and $E(G)$, respectively. 
The {\em order\/} of $G$ means the number of vetices of $G$. 
A cycle of length $k$ is a {\em $k$-cycle\/}. 
For a cycle $C$, an edge $e \in E(G)$ such that 
$V(e) \subset V(C)$ and $e \notin E(C)$ is called a {\em chord\/} of $C$. 
A cycle $C$ in $G$ is {\em separating\/} 
if $G-V(C)$ is a disconnected graph. 
We denote the induced subgraph of a graph $G$ by $S \subset V(G)$ by $G[S]$. 
A set  $M$ of edges of a graph $G$ is a {\em matching\/}
if no two edges of $M$ share a vertex. 
Let $M$ be a matching of a graph $G$.  
Each vertex incident with 
an edge of $M$ is {\em covered\/} by $M$. 
The set of vertices covered by $M$ is denoted by $V(M)$. 
In particular, 
$M$ is {\em perfect\/} if $M$ covers all vertices of $G$; that is, $V(G)=V(M)$.   
A matching $M$ of $G$ is {\em extendable\/} 
if $G$ has a perfect matching containing $M$. 
Moreover, 
a graph $G$ with at least $2k+2$ vertices is {\em $k$-extendable\/} 
if any matching $M$ in $G$ with $|M| = k$ is extendable. 
Matching extendability has been widely studied in literature (e.g., see \cite{P08}).
In particular, matching extendability of graphs on closed surfaces was investigated in 
\cite{AKP, AP11, KNPS,P1988}; for example, 
it was proven as a basic result that no planar graph is $3$-extendable.

A graph $G$ is {\em $1$-embeddable\/} on a closed surface $F^2$
if it can be drawn on $F^2$ 
so that every edge of $G$ crosses another edge at most once. 
The drawn image of $G$ on $F^2$ is a $1$-{\em embedded} graph on $F^2$. 
(We implicitly consider {\em good drawings\/}, that is, (i) vertices are on different 
points on the surface, (ii) no adjacent edges cross, (iii) no three edges cross 
at the same point, and (iv) any non-adjacent edges do not touch tangently.)
The study of $1$-{\em planar\/} graphs, which are 
$1$-embeddable graphs on the plane or the sphere, 
was first introduced by Ringel \cite{Ringel},
and recently developed in various points of view (see e.g., \cite{KLM, Book}); 
the drawn image is called a $1$-{\em plane graph\/}. 
%In particular, $1$-embeddable graph $G$ on the sphere is called {\em $1$-planar\/}. 
It is known that if $G$ is a $1$-embedded graph on $F^2$ 
with at least three vertices,
then $|E(G)| \leq 4|V(G)| - 4\chi(F^2)$ holds, 
where $\chi(F^2)$ stands for the Euler characteristic of $F^2$ 
(see \cite{Nagasawa} for example).
In particular, a $1$-embedded graph $G$ on $F^2$ that 
satisfies the equality, that is $|E(G)| = 4|V(G)| - 4\chi(F^2)$, 
is {\em optimal\/}. 
An edge in a $1$-embedded graph $G$ is {\em crossing\/} 
if it crosses another edge, 
and {\em non-crossing\/} otherwise. 
Let $G$ be an optimal $1$-embedded graph on $F^2$, and let $W$ 
be a closed walk consisting of only non-crossing edges 
that bounds a $2$-cell $D$; where a $2$-cell is homeomorphic to 
an open disc.  
If $D$ contains an odd number of vertices, then 
we call $D$ an {\em odd weighted region\/}. 
In particular, if $W$ is a cycle, then $W$ is a {\em barrier cycle\/}. 
A barrier cycle of length $k$ is called a {\em barrier $k$-cycle\/}. 

The matching extendability of $1$-embedded graphs on $F^2$ was first addressed  
in \cite{FSS}, and the authors proved that 
every optimal $1$-plane graph (abbreviated as O1PG)  of even order is $1$-extendable. 
Further in the same paper, they discussed $2$-extendability of O1PG's, and established the following theorem.

\begin{thm}[Fujisawa et al. \cite{FSS}]\label{pl:butterfly}
An O1PG $G$ of even order is $2$-extendable
unless $G$ contains a barrier $4$-cycle.
\end{thm}

Furthermore, they discussed extendable three edges in O1PG's and obtained the following result. 

\begin{thm}[Fujisawa et al. \cite{FSS}]\label{pl:5-conn}
Let $G$ be a $5$-connected O1PG of even order, 
and $M$  be a matching of $G$ with $|M| = 3$.
Then $M$ is extendable unless
$G$ contains a barrier $6$-cycle $C$ such that $V(M) = V(C)$.
\end{thm}

We extend the topic to graphs on non-spherical closed surfaces. 
In this paper, we especially discuss matching extendability of 
optimal $1$-embedded graphs on the projective plane; 
we denote the projective plane by $P^2$ briefly. 
An optimal $1$-embedded graph on $P^2$ is also 
called an {\em optimal 1-projective plane graph\/}, 
and is abbreviated as O1PPG. 
Note that every O1PPG has exactly $4|V(G)| - 4$ edges by 
the equality above with $\chi(P^2) = 1$. 
First, we discuss $1$-extendability of O1PPG's, and 
show the following theorem, using Hamiltonian paths 
contained in those graphs. 

\begin{thm}\label{thm:1-ext}
Every O1PPG of even order 
is $1$-extendable. 
\end{thm}

Next, we discuss $2$-extendability, and prove the following theorem. 
The statement looks similar to the spherical case, 
but we need to establish some lemmas specific to the case of 
the projective plane.

\begin{thm}\label{thm:2-ext}
An O1PPG $G$ of even order 
is $2$-extendable if and only if 
$G$ contains a barrier $4$-cycle.
%cycle C of length 4 s.t. G-C has an odd component
\end{thm}

The following corollary easily follows from Theorem \ref{thm:2-ext}.

\begin{cor}\label{cor:2-ext}
Any $5$-connected O1PPG of even order is $2$-extendable.
\end{cor}

Every O1PPG $G$ has a vertex with degree $6$, 
since the average degree of $G$ is less than $8$.  
(And since the minimum degree of $G$ is at least $6$, 
and every optimal $1$-embedded graph on a closed surface 
is Eulerian. We mention these facts in Section~\ref{sect:pre}.)
Therefore, no O1PPG is $3$-extendable; take three edges on the $6$-cycle 
induced by neighbors of a vertex of degree $6$.  
The following theorem characterizes three mutually independent edges that are not extendable in those graphs. 

\begin{figure}[t]
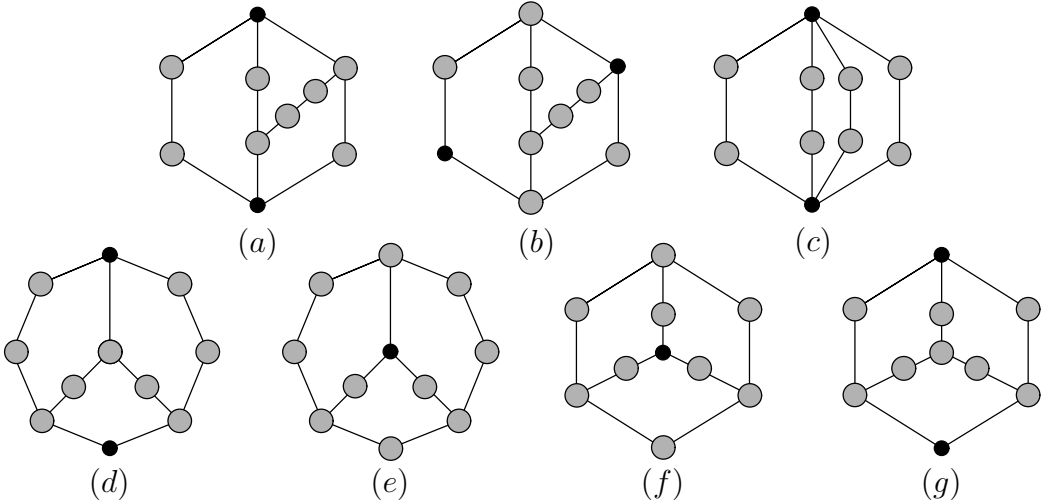

\begin{center}
{\unitlength 0.1in%
}%
\end{center}
\caption{Specified subgraphs in $Q(G)$ in Theorem~\ref{thm:3-ext}.}
\label{fig:3-ext}
\end{figure}

\begin{thm}\label{thm:3-ext}
Let $G$ be a $5$-connected O1PPG of even order, and $M$  be a matching of $G$ with $|M| = 3$.
Then $M$ is not extendable if and only if
$G$ has either 
(i) an odd weighted region bounded by a closed walk $W$ 
of length $6$ such that $V(W) \setminus V(M) = \emptyset$, or 
(ii) a subgraph of $Q(G)$ shown as (a), \ldots, (f) or (g) in Figure~\ref{fig:3-ext}, each of whose face 
is an odd weighted region, where big gray vertices are covered by $M$. 
\end{thm}

%Above figure represents graphs on $P^2$. In this figure, vertices located on diagonal 
%line represent same vertices. For example Fig.$1$ (a) is a graph on $P^2$ which 
%has seven vertices and nine edges.

Note that each of $(a)$, $(b)$, $\ldots$ , $(f)$ and $(g)$ represents a graph on $P^2$. 
To obtain the projective plane $P^2$, identify each antipodal pair of 
points of the hexagon or the octagon in the figure. 
(Similarly, carry out the same identification for regular polygons or dashed-circles in 
other figures in the paper to obtain $P^2$.)
For example $(a)$ is a graph on $P^2$ that 
has seven vertices and nine edges.

This paper is organized as follows.
In the next section, 
we first define terminology used in the paper, and 
discuss the fundamental results hold for optimal $1$-embedded 
graphs on general closed surfaces. 
Next, we 
discuss connectivity of O1PPG's and separating short cycles 
in underlying quadrangulation consisting of non-crossing edges in Section~\ref{sect:Q(G)}. 
Furthermore, we characterize O1PPG's having connectivity exactly $5$ in the section;  
note that there is no O1PG (on the sphere) having connectivity exactly $5$. 
%In Section~\ref{sect:B(G,S)}, we introduce a bipartite graph $B(G,S)$ from a graph $G$ 
%with $S\subset V(G)$, which is useful when discussing matching extendability of graphs; 
%actually we consider an O1PPG $G$, and obtain $B(G,S)$ on $P^2$. 
In Section~\ref{sect:Proofs}, 
we discuss extendability of O1PPG's, and prove our main theorems.

%%%%%%%%%%%%%%%%%%%%%%%%%%%%%%%%%%%%%%%%%%%%%%%%%%%%%%%%%%%%%%%%%%%%%%%%%%%%%%%%%

\section{Preliminaries and basic results}\label{sect:pre}

A vertex set $S$ of a connected graph $G$ is a {\em cut\/} if 
$G - S$ has at least two connected components. 
A cut $S$ of $G$ is {\em minimal\/} if any proper subset of $S$ 
is not a cut of $G$. 
For a cut $S$ of $G$, if $|S| = k$, then we call $S$ a {\em $k$-cut\/} 
of $G$. 
We denote the number of connected components of 
$G-S$ for $S\subset V(G)$ by $C(G-S)$. 
In particular, the number of {\em odd components\/} (resp., {even components\/}), 
i.e., connected components having odd (resp., even) number of vertices, 
is denoted by $C_o(G-S)$ (resp., $C_e(G-S)$). 
That is, we have $C(G-S)=C_o(G-S)+C_e(G-S)$. 

Let $G$ be a graph embedded on a closed surface $F^2$. 
Then a connected component of $F^2 - G$, 
which is as a topological space, 
is a {\em face\/} of $G$, 
and we denote the face set of $G$ by $F(G)$;  
that is, ``a face'' in this paper is not necessarily 
homeomorphic to an (open) $2$-cell. 
A {\em boundary closed walk\/} $W$ of a face $f$ is a closed walk bounding $f$ 
in $G$. 
(Actually, under our definition, the boundary of a face 
might be a union of closed walks.) %; we sometimes denote like $W = \partial f$. 
A $k$-{\em gonal face\/} or simply a $k$-{\em face\/} means a face with 
boundary closed walk of length exactly $k$. 
If every face of $G$ is homeomorphic to a $2$-cell, then 
$G$ is a {\em $2$-cell embedding\/} or {\em $2$-cell embedded graph\/} on $F^2$. 
Furthermore, a {\em region\/} bounded by a closed walk might contains 
some vertices and edges in its interior in our latter argument; 
that is, a face is always a region, but the converse 
does not hold in general.  
%Note that a $2$-cell region implies an ``open'' $2$-cell region in this paper. 

A simple closed curve $\gamma$ on a closed surface $F^2$ is 
{\em trivial\/} if $\gamma$ bounds a $2$-cell on $F^2$, and 
{\em essential\/} otherwise. 
We apply these definition to cycles of graphs embedded on $F^2$, 
regarding them as simple closed curves. 
A simple closed curve $\gamma$ on a closed surface $F^2$ 
is {\em surface separating\/} if $F^2 - \gamma$ is disconnected as a topological space. 
We also apply the definition to cycles of graphs on $F^2$. 
Note that every trivial closed curve on a closed surface 
is surface separating. 
In addition, it is well-known that every surface separating simple closed curve 
on the sphere or the projective plane is trivial. 
The following proposition is known in 
topological graph theory, and is commonly used.

\begin{pro}[Nakamoto \cite{Naka}]\label{prop:parity}
Let $G$ be a graph embedded on a closed surface $F^2$ so that each face is bounded 
by a closed walk of even length. Then the length of two cycles in $G$ have the same parity 
if they are homotopic to each other on $F^2$. 
Furthermore, there is no surface separating odd cycle in $G$. 
\end{pro}

The {\em representativity\/} $r(G)$ of
a graph $G$ embedded on a non-spherical closed surface $F^2$ is 
the minimum number of crossing points of $G$ and 
$\gamma$, 
where $\gamma$ ranges over 
all essential simple closed curves on $F^2$. 
A graph $G$ embedded on $F^2$ is {\em k-representative\/} if 
$r(G) \geq k$. 
A graph $G$ embedded on a non-spherical closed surface $F^2$ is {\em polyhedral\/} 
if $G$ is $3$-connected and $3$-representative. 
In particular, a graph $G$ embedded on the sphere is polyhedral 
if $G$ is just $3$-connected. 

A {\em quadrangulation\/} (resp., {\em triangulation\/}) is a simple $2$-cell embedded 
graph on a closed surface 
such that every face is a $4$-face (resp., $3$-face). 
It was shown in \cite{Nagasawa} that 
every simple optimal $1$-embedded graph $G$ on $F^2$ is obtained from a polyhedral 
quadrangulation $H$ by adding a pair of crossing edges in each face of $H$. 
%First we refer to a known result from literature, 
%which concerns the basic structure of optimal 1-embedded graphs on closed surfaces.  
%The following lemma shows that there is an one-to-one correspondence between the set of %polyhedral quadrangulations of $F^2$ and the set of optimal $1$-embedded graphs on  $F^2$. 
%\begin{thm}[Nagasawa et al. \cite{Nagasawa}]\label{thm:quadcross}
%Let $G$ be an optimal $1$-embedded graph on a closed surface $F^2$.
%If we remove all the crossing edges of $G$,
%then the resulting graph is a polyhedral quadrangulation of $F^2$.
%\end{thm}
We call the quadrangulation $H$, which consists of all the non-crossing edges of $G$, 
the {\em quadrangular subgraph\/} of $G$, 
and denote it by $Q(G) (=H)$. 
By the property above, $\deg_G(v)=2\deg_H(v)$ for any $v\in V(G)$, 
that is, 
$G$ is Eulerian. 
For a vertex $v$ of an optimal $1$-embedded graph $G$ on $F^2$, 
the union of all the faces (with boundaries) of $Q(G)$ incident to $v$ 
forms a disc $D$ containing the unique vertex $v$. 
We call the boundary cycle of $D$ the {\em link\/} of $v$ and 
denote it by $L_G(v)$; observe that the boundary corresponds to 
a cycle since $Q(G)$ is polyhedral. %by Theorem~\ref{thm:quadcross}. 
%For an optimal $1$-embedded graph $G$ on a closed surface, 
%the quadrangulation obtained by the above lemma 
%is denoted by $Q(G)$. 
%It should be noticed that 
%$G$ is obtained by 
%adding two diagonal edges to every face of $Q(G)$. 

Let $F^2$ be a closed surface. 
An {\em arc\/} in $F^2$ is the image of a continuous map 
$\alpha : [0, 1] \rightarrow F^2$; we denote the image $\alpha([0, 1])$ by $\alpha$ for brevity, 
if there is no misunderstanding.  
The arc $\alpha$ {\em joins\/} its endpoints $\alpha(0)$ and $\alpha(1)$. 
Let $G$ be an optimal $1$-embedded graph on $F^2$ and let $H_1$ and $H_2$ be 
connected subgraphs of $G$. 
Then a subgraph $K$ of $Q(G)$ {\em separates\/} $H_1$ and $H_2$ on $F^2$ 
if $V(K) \cap V(H_i) = \emptyset$ for each $i \in \{1, 2\}$, and any arc $\alpha$ on $F^2$ that 
joins $x_1 \in V(H_1)$ and $x_2 \in V(H_2)$ has an intersection with $K$; i.e., 
$\alpha \cap K \neq \emptyset$. 
Note that $F^2 \setminus K$ is disconnected as a 
topological space. 

Let $G$ be an optimal $1$-embedded graph on $F^2$ and 
let $\mathcal{C}(G)$ be the set of all the crossing points of $G$. 
We obtain the {\em associated graph\/} $G^\times$, which is 
embedded on $F^2$, from $G$ by regarding every crossing point as a vertex of degree $4$. 
(That is, $G^\times$ has $V(G^\times) = V(G) \cup \mathcal{C}(G)$ 
and $E(G^\times) = E(Q(G)) \cup 
\{xz, yz\,|\, xy \in E(G) \setminus E(Q(G)), z \in \mathcal{C}(G) \cap xy \}$.) 
We call new vertices, which correspond to crossing points of $G$, {\em false vertices\/} 
of $G^\times$. 
Note that $G^\times$ is a triangulation of $F^2$. 
%In the following lemma, 
%note that a face of a graph embedded on $F^2$ is not necessarily homeomorphic to 
%an open $2$-cell. 

In the argument below, we often consider the induced subgraph of 
$Q(G)$ by a cut $S$ of $G$, 
which is $Q(G)[S]$ under our definition. 
However, when the underlying graph $G$ is clear, we use $Q[S]$ in place of $Q(G)[S]$, to simplify notation. 
In the following four lemmas, we assume that 
$G$ is an optimal $1$-embedded graph on $F^2$, 
and $S\subset V(G)$ is a cut of $G$.

\begin{lm}\label{lm:sep}
Let $D_1, \ldots, D_m$ $(m \geq 2)$ denote connected components of $G - S$. 
Then any two connected components $D_i$ and $D_j$ 
$(i\neq j)$ are separated by $Q[S]$. 
That is, each face of $Q[S]$ contains at most 
one connected component of $G - S$. 
\end{lm}

\begin{Proof}
Suppose to the contrary 
that $Q[S]$ does not separate connected components $D_i$ and $D_j$ $(i \neq j)$ 
of $G-S$. 
Then there exists an arc $\alpha$ on $F^2$ such that 
$\alpha \cap Q[S] = \emptyset$ and 
$\alpha$ joins $x \in V(D_i)$ and $y \in V(D_j)$. 
%Let $G^\times$ denote the associated graph of $G$. 
Since the associated graph $G^\times$ of $G$ is a triangulation of $F^2$ and 
$\alpha \cap Q[S] = \emptyset$, 
we can fix $\alpha$ so that 
$\alpha \cap G^\times \subset V(G^\times) \setminus S$. 
We may assume that $\alpha \cap D_k = \emptyset$ 
for any $k \neq i, j$; otherwise, retake the closest $D_i$ and $D_j$. 
Since $G^\times$ is a triangulation of $F^2$, 
there exists a path $P^\times$ in $G^\times$ between $x$ and $y$ along $\alpha$. 
Now assume that $P^\times$ passes through a false vertex $z$ 
corresponding to a crossing point created by a pair of crossing 
edges $v_0v_2$ and $v_1v_3$. 
If a $2$-path $v_izv_{i + 2}$ is contained in $P^\times$, 
then we replace it by $v_iv_{i + 2}$, which is a 
crossing edge of $G$, where the indices are taken modulo $4$. 
On the other hand, if a $2$-path $v_izv_{i + 1}$ is contained in 
$P^\times$, then we replace 
it by a non-crossing edge $v_iv_{i + 1}$ of $G$. 
%where the indices are taken modulo $4$. 
We do the replacement above for all false vertices contained in 
$P^\times$, 
and obtain a path $P$ in $G$ between $x$ and $y$ not containing 
any vertex in $S$, a contradiction. 
\end{Proof}

Now, we show the following lemma that 
mentions the minimum degree of $Q[S]$ for a minimal cut $S$ 
of optimal $1$-embedded graphs on closed surfaces. 

\begin{lm}\label{lm:degree}
%Let $G$ be an optimal $1$-embedded graph on $F^2$, 
%and let $S$ be a minimal cut of $G$.
If $S$ is minimal, then the minimum degree of $Q[S]$ is at least $2$.
\end{lm}

\begin{Proof}
Let $v\in S$, and let $D_1$ and $D_2$ be connected components 
of $G-S$. 
Since $S$ is minimal, $G$ has edges $vx_1$ and $vx_2$ where 
$x_i \in V(D_i)$ for each $i \in \{1,2\}$. 
Note that both $x_1$ and $x_2$ are on the link $L_G(v)$. 
Since $L_G(v)$ is a cycle, there must be two vertices $s_1, s_2 \in S \cap V(L_G(v))$ which separate $x_1$ and $x_2$ in $L_G(v)$. 
Since $\{s_1, s_2\}$ separates $x_1$ and $x_2$ in $G$ as well, 
we have $\{s_1, s_2\} \subseteq N_{Q[S]}(v)$, 
and we got our desired conclusion. 
%Then there must be two vertices $s_1, s_2 \in S$ such 
%that each of $vs_1$ and $vs_2$ is a non-crossing edge of $G$, 
%and that $s_1$ and $s_2$ separate $x_1$ and $x_2$ on $L_G(v)$. 
%Otherwise, $x_1$ and $x_2$ can be joined by a path using 
%non-crossing edges of $L_G(v)$ and crossing edges 
%inside the disc bounded by $L_G(v)$; observe that 
%there are exactly $\frac{1}{2}|L_G(v)|$ such crossing edges. 
%Thus we got our desired conclusion. 
\end{Proof}

%In particular, the number of {\em odd components\/}, 
%i.e., connected components having odd (resp., odd) number of vertices, 
%is denoted by $o(G-S)$ in our argument below. 

\begin{lm}\label{lm:Q[S]}
If $Q[S]$ has $p$ faces such that the sum of the lengths of its boundary walks is at least $2q \geq 6$, then the following inequalities hold:
%If $Q[S]$ has $p$ faces each of whose boundary walks has (the sum of) %length 
%at least $2q \geq 6$, then the following inequalities hold: 
\begin{enumerate}
\item[(i)] $|E(Q[S])| \geq 2|F(Q[S])| + (q -2)p$ 
\item[(ii)] $|S| - \chi(F^2) + (2-q)p \geq |F(Q[S])|$ 
\end{enumerate} 
%where $\chi(F^2)$ stands for the Euler characteristic of $F^2$. 
\end{lm}

\begin{Proof}
Note that a face of $Q[S]$ is not necessarily a $2$-cell, 
and does not necessarily have the unique boundary component. 
By Proposition~\ref{prop:parity}, which actually holds for general 
boundaries of faces (that is, not restricted to cycles),  
the sum of lengths of boundary walks of $Q[S]$ is even. 
Thus, we have $2|E(Q[S])| \geq 4(|F(Q[S])|-p) + 2pq$, 
and hence (i) in the theorem holds. 
Furthermore, by combining Euler's formula 
$|S| - |E(Q[S])| + |F(Q[S])| \geq \chi(F^2)$, 
we can easily obtain (ii) in the statement. 
\end{Proof}

\begin{lm}\label{lm:edge-bound}
If $|S| \leq C_o(G-S)+2k$ holds for some integer $k$, 
then we have the following: 
$$2|F(Q[S])| + 2k - \chi(F^2) \geq |E(Q[S])|$$ 
\end{lm}

\begin{Proof}
By Lemma~\ref{lm:sep}, $|F(Q[S])| \geq C_o(G-S) \geq |S| - 2k$ holds. 
Then by $|S| - |E(Q[S])| + |F(Q[S])| \geq \chi(F^2)$, 
we obtain the inequality in the statement. 
\end{Proof}

\section{Minimal cuts and subgraphs in $Q(G)$}\label{sect:Q(G)}

In this section, we describe properties of induced subgraphs by minimal 
cuts in O1PPG's. 
First of all, we show $4$-connectedness of O1PPG's as follows. 

\begin{thm}\label{thm:4-cut}
Every O1PPG $G$ is $4$-connected. 
Furthermore, if $G$ has a $4$-cut $S$,
then $Q[S]$ contains a separating trivial $4$-cycle of $G$.
\end{thm}

\begin{Proof}
In \cite{NS}, 
it was proven that 
every quadrangulation $G$ on a closed surface with $|V(G)| \geq 6$  
can be extended to a $4$-connected triangulation by adding 
a diagonal edge in every face of $G$. 
Furthermore, 
it was shown in \cite{SuzukiPG} that 
every O1PPG has at least nine vertices. 
Combining the results above, we obtain the former half of the statement of 
the theorem.

Next, we discuss the latter half of the statement. 
Let $S$ be a $4$-cut of $G$, and first 
assume that $Q[S]$ is not a $2$-cell embedding; 
note that $Q[S]$ is bipartite. 
Since $Q[S]$ has at least two faces by Lemma~\ref{lm:sep}, 
$Q[S]$ contains a cycle $C$. 
Then $|C|=4$ and $C$ does not have any chord; otherwise $Q[S]$ would have a 
trivial cycle of length $3$, a contradiction to that $Q[S]$ is bipartite. 
In this case, $Q[S]$ is just a $4$-cycle, and has exactly two faces, one of which 
is a $2$-cell and the other of which contains a cross cap. 
Each face contains exactly one connected component by Lemma~\ref{lm:sep}, 
and we have our desired trivial $4$-cycle in $Q[S]$. 

Secondly, we assume that $Q[S]$ is a $2$-cell embedding. 
By Euler's formula $4-|E(Q(S))|+|F(Q[S])|=1$, 
$|E(Q[S])|$ is either $5$ or $6$ since $|F(Q[S])|\geq 2$ by Lemma~\ref{lm:sep}. 
In the case when $|E(Q[S])|=6$, $Q[S]\cong K_4$, 
and it is known that 
$K_4$ is uniquely embedded on $P^2$ such that each 
face is bounded by a closed walk of even length, which is actually 
a cycle of length $4$. 
On the other hand, if $|E(Q[S])|=5$, then 
$Q[S]$ has exactly two faces, one of which is bounded by a $4$-cycle 
and the other of which is bounded by a closed walk of length $6$; 
note that the sum of the lengths of boundary walks must be $2|E(Q[S])|=10$.  
Observe that in the embedding of $Q[S]$ above, any $3$-cycle must 
be essential. (Therefore, the embedding of $Q[S]$ is uniquely determined.) 
In either case, we have our desired trivial $4$-cycle in $Q[S]$ by Lemma~\ref{lm:sep}. 
%note that $G-S$ has at least two connected components.  
\end{Proof}

%According to the result in , 
%an O1PPG has at least nine vertices. 
%The following theorem implies that every O1PPG is $4$-connected. 
%
%\begin{thm}[Noguchi and Suzuki \cite{noguchi2015relationship}]
%\label{lem:4-conntri}
%Let $G$ be a quadrangulation on a closed surface $F^2$ 
%with $|V(G)| \geq 6$. 
%Then $G$ can be extended to a $4$-connected triangulation by adding 
%a diagonal edge in every face of $G$. 
%\end{thm}

Next, we present some facts holding for $5$-connected O1PPG's.

\begin{lm}\label{lm:ineq}
Let $G$ be a $5$-connected O1PPG, and let $S$ be a minimal cut with $|S|\in \{5, 6\}$. 
Then the followings hold: 
\begin{enumerate}
\item[(i)] $|E(Q[S])| \geq  2|F(Q[S])|+2$. 
\item[(ii)] If $Q[S]$ is not a $2$-cell embedding, then $|S|=6$. 
Furthermore, $Q[S]$ is a trivial $6$-cycle. 
\item[(iii)] If $|S|=6$, then $Q[S]$ has a $2$-cell face bounded by a $6$-cycle. 
\end{enumerate}
\end{lm}

\begin{Proof}
%First, note that $Q[S]$ does not contain a face bounded by 
%a closed walk with odd length; a face is not necessarily homeomorphic to a $2$-cell.   
%If $Q[S]$ contains a face $f_0$ bounded by a closed walk of length $4$, 
%then $f_0$ dose not contain a vertex of $G$, since $G$ is $5$-connected.  
%Hence $Q[S]$ has at least two faces bounded 
%by a closed walk of length at least $6$ by Lemma~\ref{lm:sep}. 
%Hence we obtain 
%\begin{equation}
%$2|E(Q[S])| \geq 4(|F(Q[S])|-2) + 6\times 2 = 4|F(Q[S])|+4$. 
%\end{equation}. 
The inequality (i) easily follows from (i) of Lemma~\ref{lm:Q[S]} with $p\geq 2$ and $q\geq 3$ 
since $C(G-S)\geq 2$, and since $G$ is $5$-connected. 
% note that every face of $Q[S]$ has a boundary walk of even length.   
Next, assume that $Q[S]$ is not a $2$-cell embedding. 
Then $Q[S]$ is a bipartite graph. 
If $|S|=5$, then $Q[S] \cong K_{2,3}$ by Lemma~\ref{lm:degree}. 
This $Q[S]$ is embedded on $P^2$ having three $4$-faces; 
observe that exactly one of them is not a $2$-cell face.  
By Lemma~\ref{lm:sep}, at least two faces above contain 
vertices of $G$, and $G$ would have a $4$-cut, a contradiction. 
Thus assume that $|S|=6$. 
In this case, either 
$Q[S] \cong K_{2,4}$ or $Q[S]$ contains a $6$-cycle. 
In the former case, there exists a $4$-face containing a vertex of $G$, 
a contradiction; similar to the case when $|S|=5$. 
Hence $Q[S]$ contains a $6$-cycle $C$. 
Under the condition, $C$ has at most two chords since $Q[S]$ is a planar graph. 
However, in any case, $Q[S]$ has at most one $k$-face with $k\geq 6$, contradicting Lemma~\ref{lm:sep}. 
Thus (ii) in the statement holds. 

Finally, we discuss (iii), and assume that $|S|=6$. 
We may assume that $Q[S]$ is a $2$-cell embedding by the result (ii) above. 
Suppose that $G$ has no $k$-face with $k=6$. 
Then $Q[S]$ has at least two faces bounded by closed walks of length at least $8$ 
by Lemma~\ref{lm:sep}.  
By (ii) of Lemma~\ref{lm:Q[S]} with $p\geq 2$ and $q\geq 4$, 
we have $1 \geq |F(Q[S])|$, contradicting $|F(Q[S])|\geq 2$.  
Thus $Q[S]$ has a $6$-face $f$ bounded by a closed walk $W$. 
Since $Q[S]$ has no chord inside $f$, $f$ contains a vertex of $G$. 
If $W$ is not a cycle, then it is contrary to $S$ being a minimal cut. 
Thus (iii) holds. 
\end{Proof}

\begin{figure}[tb]
\begin{center}
{\unitlength 0.1in%
\begin{picture}(12.5000,13.5000)(6.5000,-16.9500)%
% POLYGON 2 0 3 0 Black White  
% 8 1331 579 869 872 869 1328 1331 1621 1794 1334 1794 872 1794 872 1331 579
% 
\special{pn 8}%
\special{pa 1331 579}%
\special{pa 869 872}%
\special{pa 869 1328}%
\special{pa 1331 1621}%
\special{pa 1794 1334}%
\special{pa 1794 872}%
\special{pa 1331 579}%
\special{pa 869 872}%
\special{fp}%
% CIRCLE 2 0 0 0 Black White  
% 4 1331 1621 1331 1660 1331 1660 1331 1660
% 
\special{sh 1.000}%
\special{ia 1331 1621 39 39 0.0000000 6.2831853}%
\special{pn 8}%
\special{ar 1331 1621 39 39 0.0000000 6.2831853}%
% CIRCLE 2 0 0 0 Black White  
% 4 869 1334 869 1373 869 1373 869 1373
% 
\special{sh 1.000}%
\special{ia 869 1334 39 39 0.0000000 6.2831853}%
\special{pn 8}%
\special{ar 869 1334 39 39 0.0000000 6.2831853}%
% CIRCLE 2 0 0 0 Black White  
% 4 1794 1334 1794 1373 1794 1373 1794 1373
% 
\special{sh 1.000}%
\special{ia 1794 1334 39 39 0.0000000 6.2831853}%
\special{pn 8}%
\special{ar 1794 1334 39 39 0.0000000 6.2831853}%
% CIRCLE 2 0 0 0 Black White  
% 4 1794 872 1794 911 1794 911 1794 911
% 
\special{sh 1.000}%
\special{ia 1794 872 39 39 0.0000000 6.2831853}%
\special{pn 8}%
\special{ar 1794 872 39 39 0.0000000 6.2831853}%
% CIRCLE 2 0 0 0 Black White  
% 4 869 872 869 911 869 911 869 911
% 
\special{sh 1.000}%
\special{ia 869 872 39 39 0.0000000 6.2831853}%
\special{pn 8}%
\special{ar 869 872 39 39 0.0000000 6.2831853}%
% CIRCLE 2 0 0 0 Black White  
% 4 1331 579 1331 618 1331 618 1331 618
% 
\special{sh 1.000}%
\special{ia 1331 579 39 39 0.0000000 6.2831853}%
\special{pn 8}%
\special{ar 1331 579 39 39 0.0000000 6.2831853}%
% CIRCLE 2 0 0 0 Black White  
% 4 1330 1265 1330 1304 1330 1304 1330 1304
% 
\special{sh 1.000}%
\special{ia 1330 1265 39 39 0.0000000 6.2831853}%
\special{pn 8}%
\special{ar 1330 1265 39 39 0.0000000 6.2831853}%
% CIRCLE 2 0 0 0 Black White  
% 4 1330 935 1330 974 1330 974 1330 974
% 
\special{sh 1.000}%
\special{ia 1330 935 39 39 0.0000000 6.2831853}%
\special{pn 8}%
\special{ar 1330 935 39 39 0.0000000 6.2831853}%
% STR 2 0 3 0 Black White  
% 4 1280 390 1280 490 2 0 0 0
% $p_1$
\put(12.8000,-4.9000){\makebox(0,0)[lb]{$p_1$}}%
% STR 2 0 3 0 Black White  
% 4 1300 1740 1300 1840 2 0 0 0
% $p_1$
\put(13.0000,-18.4000){\makebox(0,0)[lb]{$p_1$}}%
% STR 2 0 3 0 Black White  
% 4 650 1330 650 1430 2 0 0 0
% $p_2$
\put(6.5000,-14.3000){\makebox(0,0)[lb]{$p_2$}}%
% STR 2 0 3 0 Black White  
% 4 1900 1330 1900 1430 2 0 0 0
% $p_3$
\put(19.0000,-14.3000){\makebox(0,0)[lb]{$p_3$}}%
% STR 2 0 3 0 Black White  
% 4 1900 730 1900 830 2 0 0 0
% $p_2$
\put(19.0000,-8.3000){\makebox(0,0)[lb]{$p_2$}}%
% STR 2 0 3 0 Black White  
% 4 700 730 700 830 2 0 0 0
% $p_3$
\put(7.0000,-8.3000){\makebox(0,0)[lb]{$p_3$}}%
% STR 2 0 3 0 Black White  
% 4 1120 1255 1120 1355 2 0 0 0
% $p_4$
\put(11.2000,-13.5500){\makebox(0,0)[lb]{$p_4$}}%
% STR 2 0 3 0 Black White  
% 4 1115 880 1115 980 2 0 0 0
% $p_5$
\put(11.1500,-9.8000){\makebox(0,0)[lb]{$p_5$}}%
% STR 2 0 3 0 Black White  
% 4 1580 995 1580 1095 5 0 0 0
% $f_1$
\put(15.8000,-10.9500){\makebox(0,0){$f_1$}}%
% LINE 2 0 3 0 Black White  
% 2 1330 575 1330 1620
% 
\special{pn 8}%
\special{pa 1330 575}%
\special{pa 1330 1620}%
\special{fp}%
\end{picture}}%
\end{center}
\caption{Projective-bowtie}
\label{fig:p-bowtie}
\end{figure}
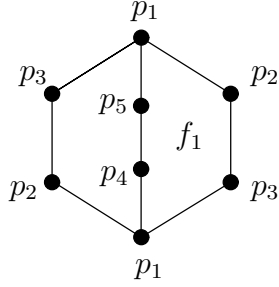

Next, see Figure~\ref{fig:p-bowtie}. 
We call the graph embedded on $P^2$ as shown in the figure a {\em projective-bowtie}.  
%A graph $G$ embedded on $P^2$ is the {\em projective-bowtie} if $G$ is isomorphic to the graph shown in Fig.$3$. 
%to obtain the projective plain $P^2$, identify each antipodal pair of points of the hexagon. 
Note that the projective-bowtie has five vertices, 
six edges and two faces bounded by closed walks of length $6$. 
%The following theorem gives a characterization of O1PPG's that have connectivity exactly $5$. 
Actually, in \cite{FSS}, 
it was proven that every $5$-connected O1PG $G$ is $6$-connected.
That is, there is no O1PG having connectivity exactly $5$. 
The next lemma (and the theorem) 
illustrates a distinct property for O1PPG's, in contrast to 
the aforementioned fact for O1PG's.

\begin{lm}\label{lm:bowtie}
Let $G$ be a $5$-connected O1PPG, and 
let $S$ be a $5$-cut of $G$. 
Then $Q[S]$ is a projective-bowtie.  
\end{lm}

\begin{Proof}
Let $S=\{p_1, \ldots, p_5\}$ be a $5$-cut of $G$; clearly, $S$ is minimal. 
By (ii) of Lemma~\ref{lm:ineq}, we may assume that 
$Q[S]$ is a $2$-cell embedding. 
It easily follows from (i) of Lemma~\ref{lm:ineq} 
that $|E(Q[S])| \geq 6$ 
by $|F(Q[S])| \ge 2$. %, $k=5$ and $C\geq 2$. 
By substituting the term of the number of faces in Euler's formula 
$5 - |E(Q[S])| + |F(Q[S])| = 1$ for (i) in Lemma~\ref{lm:ineq}, 
we have $|E(Q[S])| \leq 6$. 
Consequently, $|E(Q[S])| = 6$, and further $|F(Q[S])|=2$ holds. 
This implies that $Q[S]$ has exactly two $6$-gonal faces denoted by 
$f_1$ and $f_2$ by Lemma~\ref{lm:sep}.  
Note that $f_i$ is not bounded by a cycle for each $i\in \{1,2\}$. 
That is, there exists a vertex of $S$, say $p_1$ without loss of 
generality, which appears on the boundary closed walk 
of $f_1$, denoted by $W_1$, at least twice. 
If the distance between two $p_1$'s on $W_1$ is at most $2$, 
then either $Q[S]$ is not simple or 
$Q[S]$ has a vertex of degree $1$, contradicting 
Lemma~\ref{lm:degree}. 

Thus we may assume that 
$f_1$ is bounded by a closed walk $W=p_1xyp_1zwp_1$ where 
$\{x, y, z, w\} \subseteq \{p_2, \ldots, p_5\}$ without loss of 
generality. Observe that the $3$-cycle $p_1xyp_1$ is essential on $P^2$ 
by Proposition~\ref{prop:parity}. 
Under the condition, if $x=z$, then there exists a $2$-cell region 
bounded by a $4$-cycle $p_1wxyp_1$. 
Since $|V(Q[S])| = 5$ and $|E(Q[S])| = 6$, 
one of the two regions contains the unique inner vertex and the unique inner edge of $Q[S]$, 
contradicting Lemma~\ref{lm:degree}. 
Therefore, $x, y, z$ and $w$ are distinct vertices, and hence we 
may assume that $W=p_1p_2p_3p_1p_4p_5p_1$. 
Now, all the edges of $Q[S]$ appeared, and the outside 
region actually corresponds to the second face $f_2$ of $Q[S]$. 
That is, $Q[S]$ is the projective-bowtie. 
\end{Proof}

%Using the previous lemma, we can easily prove the following theorem. 

\begin{thm}\label{thm:6-conn}
Let $G$ be a $5$-connected O1PPG.
Then $G$ is $6$-connected if and only if $Q(G)$ does not contain 
a projective-bowtie as a subgraph. 
\end{thm}

\begin{Proof}
First, we prove the necessity. 
Assume that $Q(G)$ contains a projective-bowtie $H$ with 
$V(H) = \{p_1, p_2, p_3, p_4, p_5 \}$ as shown in 
Figure~\ref{fig:p-bowtie}, where one face of $H$, say $f_1$, 
is bounded by a closed walk 
$W = p_1p_2p_3p_1p_4p_5p_1$. 
Suppose that the region $R_1$, which corresponds to $f_1$, 
does not contain any vertex of $G$. 
Then there exists a non-crossing edge $p_2p_4$ or $p_3p_5$ of $G$ in $R_1$; now 
say $p_2p_4$, up to symmetry. 
Then there exists 
a crossing edge $p_1p_2$, which crosses $p_3p_4$ in $R_1$, contrary to $G$ being simple. 
Hence each face of $H$ contains a vertex of $G$. 
Then $V(H)$ is a $5$-cut of $G$, that is, $G$ is not $6$-connected. 
The sufficiency immediately follows from Lemma~\ref{lm:bowtie}. 
\end{Proof}

%The following lemma shows the relation 
%between a minimal $6$-cut of a $5$-connected optimal $1$-embedded graph $G$ on $P^2$
%and the structure of $Q(G)$. 

\begin{figure}[tb]
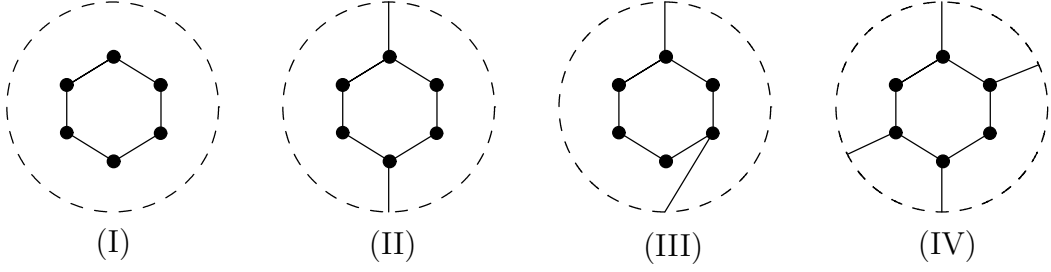

\begin{center}
{\unitlength 0.1in%
}%
\end{center}
\caption{$Q[S]$ obtained by a minimal $6$-cut $S$.}
\label{fig:6-cut}
\end{figure}

The following lemma describes the shape of $Q[S]$ for 
a minimal $6$-cut $S$ of a $5$-connected O1PPG. 
Note that in Figure~\ref{fig:6-cut}, 
antipodal points of each dashed disk are identified to 
obtain $P^2$, as we mentioned in the introduction.

\begin{lm}
\label{lm:6-cut}
Let $G$ be a $5$-connected O1PPG and $S$ be a minimal $6$-cut of $G$. Then $Q[S]$ 
is one of (I), (II), (III) and (IV) as shown in Figure~\ref{fig:6-cut}. 
\end{lm}

\begin{Proof}
By (iii) in Lemma~\ref{lm:ineq}, 
$Q[S]$ has a $2$-cell face $f$ bounded 
by a cycle $C$ of length $6$. 
By categorizing based on the number of cords of $C$ outside $f$, 
we obtain (I), (II), (III) and (IV) in Figure~\ref{fig:6-cut}. 
Note that by Lemma~\ref{lm:sep}, $Q[S]$ has at least two faces 
bounded by closed walks with length at least $6$. 
Further, consider Proposition~\ref{prop:parity}. 
\end{Proof}

\section{Proof of the main theorems}\label{sect:Proofs}

The following famous result plays an important role in the proof of 
Theorem~\ref{thm:1-ext}.

\begin{thm}[Kawarabayashi and Ozeki \cite{K-O}]\label{thm:K-O}
Every $4$-connected graph embedded on $P^2$ is Hamilton-connected. 
\end{thm}

Now, we prove our first main result in the paper. 

\bigskip

\begin{Proofof}{Theorem \ref{thm:1-ext}}
Let $G$ be an O1PPG, and let $e=uv$ be an edge of $G$. 
By the result in \cite{NS} as mentioned in the proof of Theorem~\ref{thm:4-cut}, 
$G$ has a $4$-connected triangulation $T$ 
as a spanning subgraph. 
Then $T$ has a Hamiltonian path $P=x_1x_2\cdots x_{|V(G)|}$, 
where $u=x_1$ and $v={x_{|V(G)|}}$ by Theorem~\ref{thm:K-O}. 
In $G$, $\{x_2x_3, x_4x_5, \ldots, x_{|V(G)|-2}x_{|V(G)|-1}, x_{|V(G)|}x_1\}$ is a 
perfect matching that contains $e$. 
\end{Proofof}

%Next, we show the following result that is often used in matching theory. 
Actually, the following lemma is a generalization of 
Lemma~2.3 in \cite{P1} that is often used in matching theory; 
in particular, if $k = 1$, 
then $S$ in the lemma is called an $\{e_1, e_2\}$-{\em blocker\/}.

\begin{lm}[Fujisawa et al. \cite{FSS}]\label{lem:blocker}
Let $G$ be a $k$-extendable graph and 
$\{e_1, \ldots , e_{k+1}\}$ 
be a matching of $G$ which is not extendable. 
Then there exists $S \subseteq V(G)$ 
such that 
\begin{itemize}
\item[{(i)}]
$\displaystyle S \supset \bigcup_{i=1}^{k+1} V(e_i)$ and
\item[{(ii)}]	$|S| = C_o(G - S) + 2k$. 
\end{itemize}
%where $o(H)$ stands for the number of odd components of $H$. 
\end{lm}

In the remaining part of the section, 
we prove the following two main results using tools proven in 
the previous section. 

\bigskip

\begin{Proofof}{Theorem \ref{thm:2-ext}}
The necessity is trivial and hence we prove the sufficiency of the statement. 
Let $G$ be an O1PPG that is not $2$-extendable, and assume that 
$e_1$ and $e_2$ are independent edges of $G$ that are not extendable.  
By Theorem~\ref{thm:1-ext}, 
$G$ is $1$-extendable, and hence there exists $S \subset V(G)$ 
that satisfies (i) and (ii) of Lemma~\ref{lem:blocker} for $k=1$; for $e_1$ and $e_2$.  
Now we consider $Q[S]$ on $P^2$. 
%Since every face of $Q[S]$, which is not necessarily a 
%$2$-cell, is bounded by a closed walk of even length, 
%we have $2|E(Q[S])| \geq 4(|F(Q[S])|$; that is, 
By (i) of Lemma~\ref{lm:Q[S]} with $p=0$, 
we have $|E(Q[S])| \geq 2|F(Q[S])|$. 
On the other hand, by Lemma~\ref{lm:edge-bound} with $k=1$, 
we have $2|F(Q[S])|+1 \geq |E(Q[S])|$. 
Thus, either $|E(Q[S])|=2|F(Q[S])|$ or $|E(Q[S])|=2|F(Q[S])|+1$ holds. 

In the former case, every face of $Q[S]$ is bounded by a $4$-cycle, 
and at most one of them contains a cross cap. 
That is, the other faces are all $2$-cell, 
and by Lemma~\ref{lm:sep}, we can find an odd weighted face 
since $G-S$ has at least two odd components.
In the latter case, the equality of Euler's formula (in Lemma~\ref{lm:edge-bound}) holds, 
and hence $Q[S]$ is a $2$-cell embedding. 
Furthermore, by the argument above, $Q[S]$ has the unique 
$6$-gonal face and all the others are 4-gonal.
Similarly, we find our desired barrier cycle. 
\end{Proofof}

\begin{Proofof}{Theorem~\ref{thm:3-ext}}
First, we show the necessity. 
If a $5$-connected O1PPG $G$ has (i) in the statement, 
then $G-V(M)$ has an odd component, where $M$ is the set of specified 
three independent edges. 
On the other hand, if $G$ has (ii) in the statement, 
then $G-V(M)$ has a cut vertex $v$ 
such that $G-(V(M)\cup\{v\})$ has exactly three odd components. 
In either case, $G-V(M)$ does not have a perfect matching. 

Next, we discuss the sufficiency. 
Assume that 
$M$ is a matching of a $5$-connected O1PPG $G$ with $|M|=3$ 
that is not extendable. 
By Corollary \ref{cor:2-ext}, $G$ is $2$-extendable. 
Hence there exists $S \subset V(G)$ which satisfies 
(i) and (ii) of Lemma~\ref{lem:blocker} for $k = 2$. 
Note that (i) of Lemma~\ref{lem:blocker} implies $|S|\geq 6$.

First, we show that $|S|\leq 7$. 
Suppose to the contrary that $|S| \geq 8$. 
Then $C_o(G-S) \geq 4$ by (ii) of Lemma~\ref{lem:blocker}, and hence 
%Note that since $G$ is $5$-connected, a face of $Q[S]$ 
%bounded by a closed walk with length $4$ does not contain 
%any vertex of $G$. 
$Q[S]$ has at least four faces bounded by a closed walk 
with length at least $6$ by Lemma~\ref{lm:sep}. 
Hence we obtain $|E(Q[S])| \geq 2|F(Q[S])| + 4$ by (i) of Lemma~\ref{lm:Q[S]} with 
$p\geq 4$ and $q\geq 3$. 
On the other hand, we have $2|F(Q[S])|+3 \geq |E(Q[S])|$ by Lemma~\ref{lm:edge-bound} 
and (ii) of Lemma~\ref{lem:blocker}, a contradiction. 
Therefore, we have $|S|\in \{6,7\}$. 
We will now divide the proof into the following two cases.
%In what follows, we divide the proof into the following two cases.  
%which are Case ($\alpha$) and Case ($\beta$). 

\smallskip
\noindent
Case ($\alpha$) : There exists a connected component 
$D$ of $G-S$ having exactly five neighbors in $S$. 
Let $S' \subset S$ denote the set 
of five neighbors of $D$. 
Then $Q[S']$ is the projective-bowtie by 
Lemma~\ref{lm:bowtie}.  
If $S' \subset V(M)$, then $G$ contains a region that satisfies (i) in the statement. 
Thus we assume that $S' \setminus V(M) \neq \emptyset$. 
Under the situation, note that we have $|S| = 7$; we have $C_o(G-S) = 3$ by 
Lemma~\ref{lem:blocker}. 
Then, by (ii) of Lemma~\ref{lm:Q[S]} with $q=3$ and $p=3$, 
we have $|F(Q[S])| \leq 3$, and hence $|F(Q[S])| = 3$. 
This means that the equality of (ii) of Lemma~\ref{lm:Q[S]} holds, 
and hence $Q[S]$ has exactly three $6$-gonal faces. 
Note that $|E(Q[S])| = 9$ by (i) of Lemma~\ref{lm:Q[S]} and 
Lemma~\ref{lm:edge-bound}, and that 
$C_e(G-S)=0$ by Lemma~\ref{lm:sep}.  

Now we put $S \setminus S'=\{x, y\}$. %and 
%discuss the neighbors of these vertices. 
Observe that $Q[S']$, which is a projective-bowtie, 
has exactly six edges.  
By Lemma~\ref{lm:degree}, and since $|E(Q[S])\setminus E(Q[S'])|=3$, 
we have $\deg_{Q[S]}(x)=\deg_{Q[S]}(y)=2$, and 
this implies that $s_1xys_2$ is a path of length $3$, 
where $s_1, s_2 \in S'$. 
%If $x$ and $y$ are not adjacent, that is, 
%either $x$ or $y$ is adjacent to at least two vertices in $S'$, then 
%one odd component would be contained in a region 
%bounded by a $4$-cycle, a contradiction. 
%Thus $x$ and $y$ are adjacent, and each of $x$ and $y$ has exactly one neighbor 
%in $S'$ in $Q[S]$. 
Since $Q[S]$ has exactly three $6$-gonal faces by the argument above, 
$Q[S]$ is a graph shown as (1) or (2) in Figure~\ref{fig:63region}, 
up to homeomorphism. 
Finally, we consider one vertex $z \in S \setminus V(M)$.  
If $\deg_{Q[S]}(z) = 2$, then $Q[S]$ has a region that satisfies (i) 
in the statement. 
Thus $\deg_{Q[S]}(z) \geq 3$, and hence $Q[S]$ with big gray vertices 
of $V(M)$ is (a), (b) or (c) in Figure~\ref{fig:3-ext}.

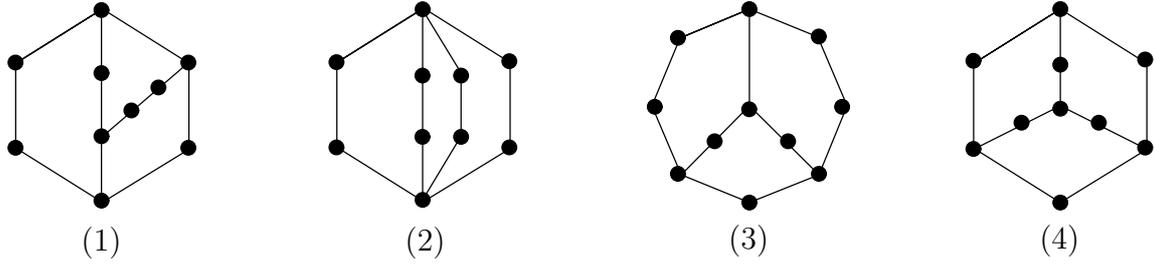
\begin{figure}[t]
\begin{center}
{\unitlength 0.1in%
\begin{picture}(59.2500,11.8500)(10.6000,-23.0500)%
% STR 2 0 3 0 Black White  
% 4 1542 2319 1542 2378 5 0 0 0
% $(1)$
\put(15.4200,-23.7800){\makebox(0,0){$(1)$}}%
% STR 2 0 3 0 Black White  
% 4 3218 2315 3218 2374 5 0 0 0
% $(2)$
\put(32.1800,-23.7400){\makebox(0,0){$(2)$}}%
% POLYGON 2 0 3 0 Black White  
% 10 4898 1158 4539 1309 4393 1663 4545 2023 4898 2168 5258 2023 5404 1663 5252 1303 5252 1303 4898 1158
% 
\special{pn 8}%
\special{pa 4898 1158}%
\special{pa 4539 1309}%
\special{pa 4393 1663}%
\special{pa 4545 2023}%
\special{pa 4898 2168}%
\special{pa 5258 2023}%
\special{pa 5404 1663}%
\special{pa 5252 1303}%
\special{pa 4898 1158}%
\special{pa 4539 1309}%
\special{fp}%
% LINE 2 0 3 0 Black White  
% 6 4898 1158 4898 1663 4545 2023 4898 1663 4898 1663 5258 2023
% 
\special{pn 8}%
\special{pa 4898 1158}%
\special{pa 4898 1663}%
\special{fp}%
\special{pa 4545 2023}%
\special{pa 4898 1663}%
\special{fp}%
\special{pa 4898 1663}%
\special{pa 5258 2023}%
\special{fp}%
% CIRCLE 2 0 0 0 Black White  
% 4 4898 1158 4898 1196 4898 1196 4898 1196
% 
\special{sh 1.000}%
\special{ia 4898 1158 38 38 0.0000000 6.2831853}%
\special{pn 8}%
\special{ar 4898 1158 38 38 0.0000000 6.2831853}%
% CIRCLE 2 0 0 0 Black White  
% 4 4898 2168 4898 2206 4898 2206 4898 2206
% 
\special{sh 1.000}%
\special{ia 4898 2168 38 38 0.0000000 6.2831853}%
\special{pn 8}%
\special{ar 4898 2168 38 38 0.0000000 6.2831853}%
% POLYGON 2 0 3 0 Black White  
% 8 6508 1158 6059 1442 6059 1884 6508 2168 6956 1890 6956 1442 6956 1442 6508 1158
% 
\special{pn 8}%
\special{pa 6508 1158}%
\special{pa 6059 1442}%
\special{pa 6059 1884}%
\special{pa 6508 2168}%
\special{pa 6956 1890}%
\special{pa 6956 1442}%
\special{pa 6508 1158}%
\special{pa 6059 1442}%
\special{fp}%
% CIRCLE 2 0 0 0 Black White  
% 4 6508 2168 6508 2206 6508 2206 6508 2206
% 
\special{sh 1.000}%
\special{ia 6508 2168 38 38 0.0000000 6.2831853}%
\special{pn 8}%
\special{ar 6508 2168 38 38 0.0000000 6.2831853}%
% CIRCLE 2 0 0 0 Black White  
% 4 6508 1158 6508 1196 6508 1196 6508 1196
% 
\special{sh 1.000}%
\special{ia 6508 1158 38 38 0.0000000 6.2831853}%
\special{pn 8}%
\special{ar 6508 1158 38 38 0.0000000 6.2831853}%
% LINE 2 0 3 0 Black White  
% 6 6508 1158 6508 1663 6508 1663 6059 1890 6508 1663 6956 1890
% 
\special{pn 8}%
\special{pa 6508 1158}%
\special{pa 6508 1663}%
\special{fp}%
\special{pa 6508 1663}%
\special{pa 6059 1890}%
\special{fp}%
\special{pa 6508 1663}%
\special{pa 6956 1890}%
\special{fp}%
% STR 2 0 3 0 Black White  
% 4 4894 2308 4894 2368 5 0 0 0
% $(3)$
\put(48.9400,-23.6800){\makebox(0,0){$(3)$}}%
% STR 2 0 3 0 Black White  
% 4 6504 2308 6504 2368 5 0 0 0
% $(4)$
\put(65.0400,-23.6800){\makebox(0,0){$(4)$}}%
% CIRCLE 2 0 0 0 Black White  
% 4 4897 1681 4897 1719 4897 1719 4897 1719
% 
\special{sh 1.000}%
\special{ia 4897 1681 38 38 0.0000000 6.2831853}%
\special{pn 8}%
\special{ar 4897 1681 38 38 0.0000000 6.2831853}%
% CIRCLE 2 0 0 0 Black White  
% 4 4528 1308 4528 1346 4528 1346 4528 1346
% 
\special{sh 1.000}%
\special{ia 4528 1308 38 38 0.0000000 6.2831853}%
\special{pn 8}%
\special{ar 4528 1308 38 38 0.0000000 6.2831853}%
% CIRCLE 2 0 0 0 Black White  
% 4 4528 2018 4528 2056 4528 2056 4528 2056
% 
\special{sh 1.000}%
\special{ia 4528 2018 38 38 0.0000000 6.2831853}%
\special{pn 8}%
\special{ar 4528 2018 38 38 0.0000000 6.2831853}%
% CIRCLE 2 0 0 0 Black White  
% 4 5258 2018 5258 2056 5258 2056 5258 2056
% 
\special{sh 1.000}%
\special{ia 5258 2018 38 38 0.0000000 6.2831853}%
\special{pn 8}%
\special{ar 5258 2018 38 38 0.0000000 6.2831853}%
% CIRCLE 2 0 0 0 Black White  
% 4 5257 1301 5257 1339 5257 1339 5257 1339
% 
\special{sh 1.000}%
\special{ia 5257 1301 38 38 0.0000000 6.2831853}%
\special{pn 8}%
\special{ar 5257 1301 38 38 0.0000000 6.2831853}%
% CIRCLE 2 0 0 0 Black White  
% 4 5378 1668 5378 1706 5378 1706 5378 1706
% 
\special{sh 1.000}%
\special{ia 5378 1668 38 38 0.0000000 6.2831853}%
\special{pn 8}%
\special{ar 5378 1668 38 38 0.0000000 6.2831853}%
% CIRCLE 2 0 0 0 Black White  
% 4 4408 1668 4408 1706 4408 1706 4408 1706
% 
\special{sh 1.000}%
\special{ia 4408 1668 38 38 0.0000000 6.2831853}%
\special{pn 8}%
\special{ar 4408 1668 38 38 0.0000000 6.2831853}%
% CIRCLE 2 0 0 0 Black White  
% 4 5098 1848 5098 1886 5098 1886 5098 1886
% 
\special{sh 1.000}%
\special{ia 5098 1848 38 38 0.0000000 6.2831853}%
\special{pn 8}%
\special{ar 5098 1848 38 38 0.0000000 6.2831853}%
% CIRCLE 2 0 0 0 Black White  
% 4 4718 1848 4718 1886 4718 1886 4718 1886
% 
\special{sh 1.000}%
\special{ia 4718 1848 38 38 0.0000000 6.2831853}%
\special{pn 8}%
\special{ar 4718 1848 38 38 0.0000000 6.2831853}%
% CIRCLE 2 0 0 0 Black White  
% 4 6508 1678 6508 1716 6508 1716 6508 1716
% 
\special{sh 1.000}%
\special{ia 6508 1678 38 38 0.0000000 6.2831853}%
\special{pn 8}%
\special{ar 6508 1678 38 38 0.0000000 6.2831853}%
% CIRCLE 2 0 0 0 Black White  
% 4 6058 1428 6058 1466 6058 1466 6058 1466
% 
\special{sh 1.000}%
\special{ia 6058 1428 38 38 0.0000000 6.2831853}%
\special{pn 8}%
\special{ar 6058 1428 38 38 0.0000000 6.2831853}%
% CIRCLE 2 0 0 0 Black White  
% 4 6508 1448 6508 1486 6508 1486 6508 1486
% 
\special{sh 1.000}%
\special{ia 6508 1448 38 38 0.0000000 6.2831853}%
\special{pn 8}%
\special{ar 6508 1448 38 38 0.0000000 6.2831853}%
% CIRCLE 2 0 0 0 Black White  
% 4 6058 1888 6058 1926 6058 1926 6058 1926
% 
\special{sh 1.000}%
\special{ia 6058 1888 38 38 0.0000000 6.2831853}%
\special{pn 8}%
\special{ar 6058 1888 38 38 0.0000000 6.2831853}%
% CIRCLE 2 0 0 0 Black White  
% 4 6307 1751 6307 1789 6307 1789 6307 1789
% 
\special{sh 1.000}%
\special{ia 6307 1751 38 38 0.0000000 6.2831853}%
\special{pn 8}%
\special{ar 6307 1751 38 38 0.0000000 6.2831853}%
% CIRCLE 2 0 0 0 Black White  
% 4 6707 1751 6707 1789 6707 1789 6707 1789
% 
\special{sh 1.000}%
\special{ia 6707 1751 38 38 0.0000000 6.2831853}%
\special{pn 8}%
\special{ar 6707 1751 38 38 0.0000000 6.2831853}%
% CIRCLE 2 0 0 0 Black White  
% 4 6947 1421 6947 1459 6947 1459 6947 1459
% 
\special{sh 1.000}%
\special{ia 6947 1421 38 38 0.0000000 6.2831853}%
\special{pn 8}%
\special{ar 6947 1421 38 38 0.0000000 6.2831853}%
% CIRCLE 2 0 0 0 Black White  
% 4 6947 1881 6947 1919 6947 1919 6947 1919
% 
\special{sh 1.000}%
\special{ia 6947 1881 38 38 0.0000000 6.2831853}%
\special{pn 8}%
\special{ar 6947 1881 38 38 0.0000000 6.2831853}%
% LINE 2 0 3 0 Black White  
% 2 1543 1828 1993 1438
% 
\special{pn 8}%
\special{pa 1543 1828}%
\special{pa 1993 1438}%
\special{fp}%
% LINE 2 0 3 0 Black Black  
% 2 1543 1153 1543 2168
% 
\special{pn 8}%
\special{pa 1543 1153}%
\special{pa 1543 2168}%
\special{fp}%
% POLYGON 2 0 3 0 Black Black  
% 8 1548 1157 1099 1441 1099 1883 1548 2167 1996 1890 1996 1441 1996 1441 1548 1157
% 
\special{pn 8}%
\special{pa 1548 1157}%
\special{pa 1099 1441}%
\special{pa 1099 1883}%
\special{pa 1548 2167}%
\special{pa 1996 1890}%
\special{pa 1996 1441}%
\special{pa 1548 1157}%
\special{pa 1099 1441}%
\special{fp}%
% CIRCLE 2 0 0 0 Black Black  
% 4 1543 1163 1543 1201 1543 1201 1543 1201
% 
\special{sh 1.000}%
\special{ia 1543 1163 38 38 0.0000000 6.2831853}%
\special{pn 8}%
\special{ar 1543 1163 38 38 0.0000000 6.2831853}%
% CIRCLE 2 0 0 0 Black Black  
% 4 1543 2158 1543 2196 1543 2196 1543 2196
% 
\special{sh 1.000}%
\special{ia 1543 2158 38 38 0.0000000 6.2831853}%
\special{pn 8}%
\special{ar 1543 2158 38 38 0.0000000 6.2831853}%
% LINE 2 0 3 0 Black Black  
% 2 3206 1149 3206 2164
% 
\special{pn 8}%
\special{pa 3206 1149}%
\special{pa 3206 2164}%
\special{fp}%
% POLYGON 2 0 3 0 Black Black  
% 8 3211 1153 2762 1437 2762 1879 3211 2163 3659 1886 3659 1437 3659 1437 3211 1153
% 
\special{pn 8}%
\special{pa 3211 1153}%
\special{pa 2762 1437}%
\special{pa 2762 1879}%
\special{pa 3211 2163}%
\special{pa 3659 1886}%
\special{pa 3659 1437}%
\special{pa 3211 1153}%
\special{pa 2762 1437}%
\special{fp}%
% CIRCLE 2 0 0 0 Black Black  
% 4 3206 1159 3206 1197 3206 1197 3206 1197
% 
\special{sh 1.000}%
\special{ia 3206 1159 38 38 0.0000000 6.2831853}%
\special{pn 8}%
\special{ar 3206 1159 38 38 0.0000000 6.2831853}%
% CIRCLE 2 0 0 0 Black Black  
% 4 3206 2154 3206 2192 3206 2192 3206 2192
% 
\special{sh 1.000}%
\special{ia 3206 2154 38 38 0.0000000 6.2831853}%
\special{pn 8}%
\special{ar 3206 2154 38 38 0.0000000 6.2831853}%
% LINE 2 0 3 0 Black White  
% 6 3205 1161 3405 1491 3405 1491 3405 1821 3405 1821 3205 2156
% 
\special{pn 8}%
\special{pa 3205 1161}%
\special{pa 3405 1491}%
\special{fp}%
\special{pa 3405 1491}%
\special{pa 3405 1821}%
\special{fp}%
\special{pa 3405 1821}%
\special{pa 3205 2156}%
\special{fp}%
% CIRCLE 2 0 0 0 Black Black  
% 4 1098 1437 1098 1475 1098 1475 1098 1475
% 
\special{sh 1.000}%
\special{ia 1098 1437 38 38 0.0000000 6.2831853}%
\special{pn 8}%
\special{ar 1098 1437 38 38 0.0000000 6.2831853}%
% CIRCLE 2 0 0 0 Black Black  
% 4 1098 1882 1098 1920 1098 1920 1098 1920
% 
\special{sh 1.000}%
\special{ia 1098 1882 38 38 0.0000000 6.2831853}%
\special{pn 8}%
\special{ar 1098 1882 38 38 0.0000000 6.2831853}%
% CIRCLE 2 0 0 0 Black Black  
% 4 1993 1882 1993 1920 1993 1920 1993 1920
% 
\special{sh 1.000}%
\special{ia 1993 1882 38 38 0.0000000 6.2831853}%
\special{pn 8}%
\special{ar 1993 1882 38 38 0.0000000 6.2831853}%
% CIRCLE 2 0 0 0 Black Black  
% 4 1993 1437 1993 1475 1993 1475 1993 1475
% 
\special{sh 1.000}%
\special{ia 1993 1437 38 38 0.0000000 6.2831853}%
\special{pn 8}%
\special{ar 1993 1437 38 38 0.0000000 6.2831853}%
% CIRCLE 2 0 0 0 Black Black  
% 4 1543 1822 1543 1860 1543 1860 1543 1860
% 
\special{sh 1.000}%
\special{ia 1543 1822 38 38 0.0000000 6.2831853}%
\special{pn 8}%
\special{ar 1543 1822 38 38 0.0000000 6.2831853}%
% CIRCLE 2 0 0 0 Black Black  
% 4 1543 1492 1543 1530 1543 1530 1543 1530
% 
\special{sh 1.000}%
\special{ia 1543 1492 38 38 0.0000000 6.2831853}%
\special{pn 8}%
\special{ar 1543 1492 38 38 0.0000000 6.2831853}%
% CIRCLE 2 0 0 0 Black Black  
% 4 1698 1687 1698 1725 1698 1725 1698 1725
% 
\special{sh 1.000}%
\special{ia 1698 1687 38 38 0.0000000 6.2831853}%
\special{pn 8}%
\special{ar 1698 1687 38 38 0.0000000 6.2831853}%
% CIRCLE 2 0 0 0 Black Black  
% 4 1838 1567 1838 1605 1838 1605 1838 1605
% 
\special{sh 1.000}%
\special{ia 1838 1567 38 38 0.0000000 6.2831853}%
\special{pn 8}%
\special{ar 1838 1567 38 38 0.0000000 6.2831853}%
% CIRCLE 2 0 0 0 Black Black  
% 4 2760 1435 2760 1473 2760 1473 2760 1473
% 
\special{sh 1.000}%
\special{ia 2760 1435 38 38 0.0000000 6.2831853}%
\special{pn 8}%
\special{ar 2760 1435 38 38 0.0000000 6.2831853}%
% CIRCLE 2 0 0 0 Black Black  
% 4 2760 1880 2760 1918 2760 1918 2760 1918
% 
\special{sh 1.000}%
\special{ia 2760 1880 38 38 0.0000000 6.2831853}%
\special{pn 8}%
\special{ar 2760 1880 38 38 0.0000000 6.2831853}%
% CIRCLE 2 0 0 0 Black Black  
% 4 3655 1430 3655 1468 3655 1468 3655 1468
% 
\special{sh 1.000}%
\special{ia 3655 1430 38 38 0.0000000 6.2831853}%
\special{pn 8}%
\special{ar 3655 1430 38 38 0.0000000 6.2831853}%
% CIRCLE 2 0 0 0 Black Black  
% 4 3655 1880 3655 1918 3655 1918 3655 1918
% 
\special{sh 1.000}%
\special{ia 3655 1880 38 38 0.0000000 6.2831853}%
\special{pn 8}%
\special{ar 3655 1880 38 38 0.0000000 6.2831853}%
% CIRCLE 2 0 0 0 Black Black  
% 4 3405 1505 3405 1543 3405 1543 3405 1543
% 
\special{sh 1.000}%
\special{ia 3405 1505 38 38 0.0000000 6.2831853}%
\special{pn 8}%
\special{ar 3405 1505 38 38 0.0000000 6.2831853}%
% CIRCLE 2 0 0 0 Black Black  
% 4 3405 1825 3405 1863 3405 1863 3405 1863
% 
\special{sh 1.000}%
\special{ia 3405 1825 38 38 0.0000000 6.2831853}%
\special{pn 8}%
\special{ar 3405 1825 38 38 0.0000000 6.2831853}%
% CIRCLE 2 0 0 0 Black Black  
% 4 3205 1825 3205 1863 3205 1863 3205 1863
% 
\special{sh 1.000}%
\special{ia 3205 1825 38 38 0.0000000 6.2831853}%
\special{pn 8}%
\special{ar 3205 1825 38 38 0.0000000 6.2831853}%
% CIRCLE 2 0 0 0 Black Black  
% 4 3205 1505 3205 1543 3205 1543 3205 1543
% 
\special{sh 1.000}%
\special{ia 3205 1505 38 38 0.0000000 6.2831853}%
\special{pn 8}%
\special{ar 3205 1505 38 38 0.0000000 6.2831853}%
\end{picture}}%
\end{center}
\caption{$Q[S]$ in the case of $|S| = 7$.}
\label{fig:63region}
\end{figure}

\smallskip
\noindent
Case ($\beta$) :  There exists no connected component 
of $G-S$ having exactly five neighbors in $S$. 
First, assume that $|S| = 6$, that is, $S = V(M)$. 
Then by Lemma~\ref{lem:blocker}, $C_o(G-S)=2$. 
By the hypothesis of Case ($\beta$), $S$ is a minimal $6$-cut of $G$. 
By Lemma~\ref{lm:6-cut}, 
$Q[S]$ is one of (I), (II), (III) and (IV) shown in Figure~\ref{fig:6-cut}.  
Then $Q[S]$ contains a region that satisfies (i) in the statement.

Next, assume that $|S| = 7$. 
By Lemma~\ref{lem:blocker}, $C_o(G-S) = 3$. 
Then by the same argument as in Case ($\alpha$),
%Next, assume that $|S| = 7$, and hence $C_o(G-S) = 3$. 
%Actually, the same argument as that with $|S|=7$ in Case ($\alpha$) holds, 
we have $|F(Q[S])|=3$ and $|E(Q[S])|=9$; observe that 
the equality in (ii) of Lemma~\ref{lm:Q[S]} holds, 
and hence $Q[S]$ is a $2$-cell embedding. 
Let $f_1, f_2$ and $f_3$ denote the three faces of $Q[S]$, 
each of which is a $6$-gonal face.  
%(Note that in this stage, we do not know whether $Q[S]$ is a $2$-cell embedding or not. 
%So, we have to use the inequality $7 - |E(Q[S])| + |F(Q[S])| = 1$ in this case. 
%However, we obtain the numbers of faces and edges above as well as the former argument. 
%Consequently, every equality holds, and we know that $Q[S]$ is a $2$-cell embedding.) 
Then, $f_i$ contains exactly one odd component of $G-S$ for each 
$i\in \{1,2,3\}$ by Lemma~\ref{lm:sep}. 
Note that $f_i$ is bounded by a $6$-cycle for 
each $i\in \{1,2,3\}$; 
otherwise, there would be an odd component satisfying the 
condition of Case ($\alpha$). 

Let $C$ denote  the boundary $6$-cycle of $f_1$, 
and put $S'=V(C)$.  
Then $Q[S']$ is (I), (II), (III) or (IV) in Figure~\ref{fig:6-cut} 
by Lemma~\ref{lm:6-cut}, up to homeomorphism. 
Here, denote the unique vertex in $S\setminus S'$ by $v$. 
If $Q[S']$ is (I) of Figure~\ref{fig:6-cut}, 
then $v$ is adjacent to exactly 
three vertices of $S'$ since $|E(Q[S])\setminus E(Q[S'])|=3$. 
Although $\deg_{Q[S]}(v)=3$, only $f_2$ and $f_3$ 
can be incident to $v$. 
That is, $v$ appears twice on the boundary closed walk of either $f_2$ or $f_3$, 
contradicting that every face of $Q[S]$ is bounded by a cycle. 
Moreover, (IV) of Figure~\ref{fig:6-cut} is not the case, 
since $v$ is adjacent to exactly one vertex of $V(C)$, contradicting 
Lemma~\ref{lm:degree}. 

Therefore $Q[S']$ is either (II) or (III) of Figure~\ref{fig:6-cut}. 
Observe that (II) is bipartite and (III) is non-bipartite. 
By Proposition~\ref{prop:parity}, if $Q[S']$ is (II) (resp., (III)), 
then $f_2$ and $f_3$ also have configurations (II) (resp., (III)). 
See Figure~\ref{fig:octa}. 
The center configuration illustrates (II) in an alternative form. 
There are exactly two ways, up to symmetry, to add $v$ inside the unique 
$8$-gonal face such that $\deg(v)=2$ and 
the $8$-gonal face is divided into two $6$-gonal faces; 
see the left-hand side one and the right-hand side one in the figure. 
However, the configuration on the right-hand side has multiple edges incident to $v$, 
a contradiction. 
Therefore, $Q[S]$ corresponds to the configuration on the left-hand side, 
which is (3) in Figure~\ref{fig:63region}. 
Similarly, (4) is derived from (III); the cases that do not lead to (4) contradict our assumptions. 
%there is the unique way to add $v$ into the $8$-gonal face of (II) or (III), 
%up to homeomorphism.   
Finally, we specify a small black vertex, which is not in $V(M)$ as well 
as Case ($\alpha$), and obtain (d), (e), (f) and (g) in Figure~\ref{fig:3-ext}. 
\end{Proofof}

\begin{figure}[t]
\begin{center}
{\unitlength 0.1in%
\begin{picture}(49.4600,11.7700)(4.7000,-14.9700)%
% POLYGON 2 0 3 0 Black White  
% 10 2898 358 2539 509 2393 863 2545 1223 2898 1368 3258 1223 3404 863 3252 503 3252 503 2898 358
% 
\special{pn 8}%
\special{pa 2898 358}%
\special{pa 2539 509}%
\special{pa 2393 863}%
\special{pa 2545 1223}%
\special{pa 2898 1368}%
\special{pa 3258 1223}%
\special{pa 3404 863}%
\special{pa 3252 503}%
\special{pa 2898 358}%
\special{pa 2539 509}%
\special{fp}%
% CIRCLE 2 0 0 0 Black White  
% 4 2898 358 2898 396 2898 396 2898 396
% 
\special{sh 1.000}%
\special{ia 2898 358 38 38 0.0000000 6.2831853}%
\special{pn 8}%
\special{ar 2898 358 38 38 0.0000000 6.2831853}%
% CIRCLE 2 0 0 0 Black White  
% 4 2898 1368 2898 1406 2898 1406 2898 1406
% 
\special{sh 1.000}%
\special{ia 2898 1368 38 38 0.0000000 6.2831853}%
\special{pn 8}%
\special{ar 2898 1368 38 38 0.0000000 6.2831853}%
% STR 2 0 3 0 Black White  
% 4 2894 1508 2894 1568 5 0 0 0
% (II)
\put(28.9400,-15.6800){\makebox(0,0){(II)}}%
% CIRCLE 2 0 0 0 Black White  
% 4 2895 655 2895 693 2895 693 2895 693
% 
\special{sh 1.000}%
\special{ia 2895 655 38 38 0.0000000 6.2831853}%
\special{pn 8}%
\special{ar 2895 655 38 38 0.0000000 6.2831853}%
% CIRCLE 2 0 0 0 Black White  
% 4 2528 508 2528 546 2528 546 2528 546
% 
\special{sh 1.000}%
\special{ia 2528 508 38 38 0.0000000 6.2831853}%
\special{pn 8}%
\special{ar 2528 508 38 38 0.0000000 6.2831853}%
% CIRCLE 2 0 0 0 Black White  
% 4 2528 1218 2528 1256 2528 1256 2528 1256
% 
\special{sh 1.000}%
\special{ia 2528 1218 38 38 0.0000000 6.2831853}%
\special{pn 8}%
\special{ar 2528 1218 38 38 0.0000000 6.2831853}%
% CIRCLE 2 0 0 0 Black White  
% 4 3258 1218 3258 1256 3258 1256 3258 1256
% 
\special{sh 1.000}%
\special{ia 3258 1218 38 38 0.0000000 6.2831853}%
\special{pn 8}%
\special{ar 3258 1218 38 38 0.0000000 6.2831853}%
% CIRCLE 2 0 0 0 Black White  
% 4 3257 501 3257 539 3257 539 3257 539
% 
\special{sh 1.000}%
\special{ia 3257 501 38 38 0.0000000 6.2831853}%
\special{pn 8}%
\special{ar 3257 501 38 38 0.0000000 6.2831853}%
% CIRCLE 2 0 0 0 Black White  
% 4 3378 868 3378 906 3378 906 3378 906
% 
\special{sh 1.000}%
\special{ia 3378 868 38 38 0.0000000 6.2831853}%
\special{pn 8}%
\special{ar 3378 868 38 38 0.0000000 6.2831853}%
% CIRCLE 2 0 0 0 Black White  
% 4 2408 868 2408 906 2408 906 2408 906
% 
\special{sh 1.000}%
\special{ia 2408 868 38 38 0.0000000 6.2831853}%
\special{pn 8}%
\special{ar 2408 868 38 38 0.0000000 6.2831853}%
% CIRCLE 2 0 0 0 Black White  
% 4 2770 1020 2770 1058 2770 1058 2770 1058
% 
\special{sh 1.000}%
\special{ia 2770 1020 38 38 0.0000000 6.2831853}%
\special{pn 8}%
\special{ar 2770 1020 38 38 0.0000000 6.2831853}%
% LINE 2 0 3 0 Black White  
% 6 2525 1215 2775 1020 2775 1020 2900 665 2900 665 2900 365
% 
\special{pn 8}%
\special{pa 2525 1215}%
\special{pa 2775 1020}%
\special{fp}%
\special{pa 2775 1020}%
\special{pa 2900 665}%
\special{fp}%
\special{pa 2900 665}%
\special{pa 2900 365}%
\special{fp}%
% STR 2 0 3 0 Black White  
% 4 2685 690 2685 750 5 0 0 0
% $f_1$
\put(26.8500,-7.5000){\makebox(0,0){$f_1$}}%
% VECTOR 2 0 3 0 Black White  
% 2 2160 870 1745 870
% 
\special{pn 8}%
\special{pa 2160 870}%
\special{pa 1745 870}%
\special{fp}%
\special{sh 1}%
\special{pa 1745 870}%
\special{pa 1812 890}%
\special{pa 1798 870}%
\special{pa 1812 850}%
\special{pa 1745 870}%
\special{fp}%
% VECTOR 2 0 3 0 Black White  
% 2 3640 870 4055 870
% 
\special{pn 8}%
\special{pa 3640 870}%
\special{pa 4055 870}%
\special{fp}%
\special{sh 1}%
\special{pa 4055 870}%
\special{pa 3988 850}%
\special{pa 4002 870}%
\special{pa 3988 890}%
\special{pa 4055 870}%
\special{fp}%
% POLYGON 2 0 3 0 Black White  
% 10 998 358 639 509 493 863 645 1223 998 1368 1358 1223 1504 863 1352 503 1352 503 998 358
% 
\special{pn 8}%
\special{pa 998 358}%
\special{pa 639 509}%
\special{pa 493 863}%
\special{pa 645 1223}%
\special{pa 998 1368}%
\special{pa 1358 1223}%
\special{pa 1504 863}%
\special{pa 1352 503}%
\special{pa 998 358}%
\special{pa 639 509}%
\special{fp}%
% CIRCLE 2 0 0 0 Black White  
% 4 998 358 998 396 998 396 998 396
% 
\special{sh 1.000}%
\special{ia 998 358 38 38 0.0000000 6.2831853}%
\special{pn 8}%
\special{ar 998 358 38 38 0.0000000 6.2831853}%
% CIRCLE 2 0 0 0 Black White  
% 4 998 1368 998 1406 998 1406 998 1406
% 
\special{sh 1.000}%
\special{ia 998 1368 38 38 0.0000000 6.2831853}%
\special{pn 8}%
\special{ar 998 1368 38 38 0.0000000 6.2831853}%
% CIRCLE 2 0 0 0 Black White  
% 4 995 655 995 693 995 693 995 693
% 
\special{sh 1.000}%
\special{ia 995 655 38 38 0.0000000 6.2831853}%
\special{pn 8}%
\special{ar 995 655 38 38 0.0000000 6.2831853}%
% CIRCLE 2 0 0 0 Black White  
% 4 628 508 628 546 628 546 628 546
% 
\special{sh 1.000}%
\special{ia 628 508 38 38 0.0000000 6.2831853}%
\special{pn 8}%
\special{ar 628 508 38 38 0.0000000 6.2831853}%
% CIRCLE 2 0 0 0 Black White  
% 4 628 1218 628 1256 628 1256 628 1256
% 
\special{sh 1.000}%
\special{ia 628 1218 38 38 0.0000000 6.2831853}%
\special{pn 8}%
\special{ar 628 1218 38 38 0.0000000 6.2831853}%
% CIRCLE 2 0 0 0 Black White  
% 4 1358 1218 1358 1256 1358 1256 1358 1256
% 
\special{sh 1.000}%
\special{ia 1358 1218 38 38 0.0000000 6.2831853}%
\special{pn 8}%
\special{ar 1358 1218 38 38 0.0000000 6.2831853}%
% CIRCLE 2 0 0 0 Black White  
% 4 1357 501 1357 539 1357 539 1357 539
% 
\special{sh 1.000}%
\special{ia 1357 501 38 38 0.0000000 6.2831853}%
\special{pn 8}%
\special{ar 1357 501 38 38 0.0000000 6.2831853}%
% CIRCLE 2 0 0 0 Black White  
% 4 1478 868 1478 906 1478 906 1478 906
% 
\special{sh 1.000}%
\special{ia 1478 868 38 38 0.0000000 6.2831853}%
\special{pn 8}%
\special{ar 1478 868 38 38 0.0000000 6.2831853}%
% CIRCLE 2 0 0 0 Black White  
% 4 508 868 508 906 508 906 508 906
% 
\special{sh 1.000}%
\special{ia 508 868 38 38 0.0000000 6.2831853}%
\special{pn 8}%
\special{ar 508 868 38 38 0.0000000 6.2831853}%
% CIRCLE 2 0 0 0 Black White  
% 4 870 1020 870 1058 870 1058 870 1058
% 
\special{sh 1.000}%
\special{ia 870 1020 38 38 0.0000000 6.2831853}%
\special{pn 8}%
\special{ar 870 1020 38 38 0.0000000 6.2831853}%
% LINE 2 0 3 0 Black White  
% 6 625 1215 875 1020 875 1020 1000 665 1000 665 1000 365
% 
\special{pn 8}%
\special{pa 625 1215}%
\special{pa 875 1020}%
\special{fp}%
\special{pa 875 1020}%
\special{pa 1000 665}%
\special{fp}%
\special{pa 1000 665}%
\special{pa 1000 365}%
\special{fp}%
% LINE 2 0 3 0 Black White  
% 2 1000 645 1350 1210
% 
\special{pn 8}%
\special{pa 1000 645}%
\special{pa 1350 1210}%
\special{fp}%
% CIRCLE 2 0 0 0 Black White  
% 4 1190 940 1190 978 1190 978 1190 978
% 
\special{sh 1.000}%
\special{ia 1190 940 38 38 0.0000000 6.2831853}%
\special{pn 8}%
\special{ar 1190 940 38 38 0.0000000 6.2831853}%
% STR 2 0 3 0 Black White  
% 4 995 1510 995 1570 5 0 0 0
% $(3)$
\put(9.9500,-15.7000){\makebox(0,0){$(3)$}}%
% POLYGON 2 0 3 0 Black White  
% 10 4898 358 4539 509 4393 863 4545 1223 4898 1368 5258 1223 5404 863 5252 503 5252 503 4898 358
% 
\special{pn 8}%
\special{pa 4898 358}%
\special{pa 4539 509}%
\special{pa 4393 863}%
\special{pa 4545 1223}%
\special{pa 4898 1368}%
\special{pa 5258 1223}%
\special{pa 5404 863}%
\special{pa 5252 503}%
\special{pa 4898 358}%
\special{pa 4539 509}%
\special{fp}%
% CIRCLE 2 0 0 0 Black White  
% 4 4898 358 4898 396 4898 396 4898 396
% 
\special{sh 1.000}%
\special{ia 4898 358 38 38 0.0000000 6.2831853}%
\special{pn 8}%
\special{ar 4898 358 38 38 0.0000000 6.2831853}%
% CIRCLE 2 0 0 0 Black White  
% 4 4898 1368 4898 1406 4898 1406 4898 1406
% 
\special{sh 1.000}%
\special{ia 4898 1368 38 38 0.0000000 6.2831853}%
\special{pn 8}%
\special{ar 4898 1368 38 38 0.0000000 6.2831853}%
% CIRCLE 2 0 0 0 Black White  
% 4 4895 655 4895 693 4895 693 4895 693
% 
\special{sh 1.000}%
\special{ia 4895 655 38 38 0.0000000 6.2831853}%
\special{pn 8}%
\special{ar 4895 655 38 38 0.0000000 6.2831853}%
% CIRCLE 2 0 0 0 Black White  
% 4 4528 508 4528 546 4528 546 4528 546
% 
\special{sh 1.000}%
\special{ia 4528 508 38 38 0.0000000 6.2831853}%
\special{pn 8}%
\special{ar 4528 508 38 38 0.0000000 6.2831853}%
% CIRCLE 2 0 0 0 Black White  
% 4 4528 1218 4528 1256 4528 1256 4528 1256
% 
\special{sh 1.000}%
\special{ia 4528 1218 38 38 0.0000000 6.2831853}%
\special{pn 8}%
\special{ar 4528 1218 38 38 0.0000000 6.2831853}%
% CIRCLE 2 0 0 0 Black White  
% 4 5258 1218 5258 1256 5258 1256 5258 1256
% 
\special{sh 1.000}%
\special{ia 5258 1218 38 38 0.0000000 6.2831853}%
\special{pn 8}%
\special{ar 5258 1218 38 38 0.0000000 6.2831853}%
% CIRCLE 2 0 0 0 Black White  
% 4 5257 501 5257 539 5257 539 5257 539
% 
\special{sh 1.000}%
\special{ia 5257 501 38 38 0.0000000 6.2831853}%
\special{pn 8}%
\special{ar 5257 501 38 38 0.0000000 6.2831853}%
% CIRCLE 2 0 0 0 Black White  
% 4 5378 868 5378 906 5378 906 5378 906
% 
\special{sh 1.000}%
\special{ia 5378 868 38 38 0.0000000 6.2831853}%
\special{pn 8}%
\special{ar 5378 868 38 38 0.0000000 6.2831853}%
% CIRCLE 2 0 0 0 Black White  
% 4 4408 868 4408 906 4408 906 4408 906
% 
\special{sh 1.000}%
\special{ia 4408 868 38 38 0.0000000 6.2831853}%
\special{pn 8}%
\special{ar 4408 868 38 38 0.0000000 6.2831853}%
% CIRCLE 2 0 0 0 Black White  
% 4 4770 1020 4770 1058 4770 1058 4770 1058
% 
\special{sh 1.000}%
\special{ia 4770 1020 38 38 0.0000000 6.2831853}%
\special{pn 8}%
\special{ar 4770 1020 38 38 0.0000000 6.2831853}%
% LINE 2 0 3 0 Black White  
% 6 4525 1215 4775 1020 4775 1020 4900 665 4900 665 4900 365
% 
\special{pn 8}%
\special{pa 4525 1215}%
\special{pa 4775 1020}%
\special{fp}%
\special{pa 4775 1020}%
\special{pa 4900 665}%
\special{fp}%
\special{pa 4900 665}%
\special{pa 4900 365}%
\special{fp}%
% LINE 2 0 3 0 Black White  
% 4 4900 365 5060 900 5060 900 4895 1365
% 
\special{pn 8}%
\special{pa 4900 365}%
\special{pa 5060 900}%
\special{fp}%
\special{pa 5060 900}%
\special{pa 4895 1365}%
\special{fp}%
% CIRCLE 2 0 0 0 Black White  
% 4 5060 895 5060 933 5060 933 5060 933
% 
\special{sh 1.000}%
\special{ia 5060 895 38 38 0.0000000 6.2831853}%
\special{pn 8}%
\special{ar 5060 895 38 38 0.0000000 6.2831853}%
% STR 2 0 3 0 Black White  
% 4 1275 810 1275 870 5 0 0 0
% $v$
\put(12.7500,-8.7000){\makebox(0,0){$v$}}%
% STR 2 0 3 0 Black White  
% 4 5145 770 5145 830 5 0 0 0
% $v$
\put(51.4500,-8.3000){\makebox(0,0){$v$}}%
% STR 2 0 3 0 Black White  
% 4 785 690 785 750 5 0 0 0
% $f_1$
\put(7.8500,-7.5000){\makebox(0,0){$f_1$}}%
% STR 2 0 3 0 Black White  
% 4 4685 690 4685 750 5 0 0 0
% $f_1$
\put(46.8500,-7.5000){\makebox(0,0){$f_1$}}%
\end{picture}}%

\end{center}
\caption{Adding $v$ of degree $2$ inside the $8$-gonal region.}
\label{fig:octa}
\end{figure}
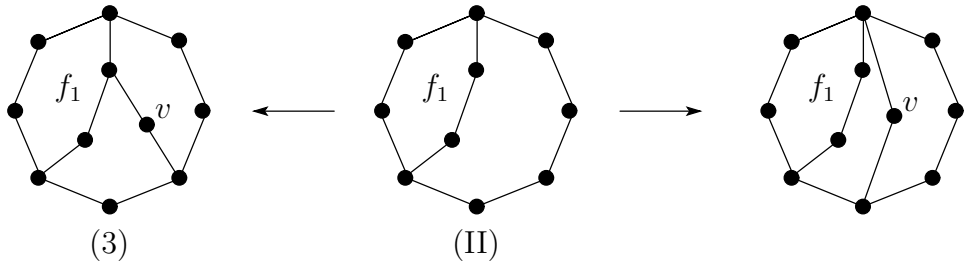

\section{Remarks}\label{sect:remarks}
In this paper, we have discussed matching extendability of 
O1PPG's. 
Then, how about for optimal $1$-embedded graphs on the 
torus (or the Klein bottle)? 
We are aware that the situation changes significantly, 
particularly in contrast to the discussion on the sphere and the projective plane. 
At least, there exist optimal 
$1$-embedded graphs on the torus that are not $1$-extendable. 
(Consider an optimal $1$-embedded graph $G$ on the torus 
such that $Q(G)$ has a subgraph $H$ that is also a quadrangulation of the torus, 
and that every face of $H$ is an odd weighted region of $G$. 
Then every non-crossing edge of $H$ is not extendable; 
observe that $|V(H)|=C_o(G-V(H))$.) 
%We know that there exists a famous conjecture by Gr\"{u}nbaum \cite{Grun} and 
%independently Nash-Williams \cite{Nash}, which asserts that every $4$-connected 
%toroidal graph has a Hamilton cycle; Kawarabayashi and Ozeki \cite{KandO} 
%solved it for $4$-connected toroidal triangulations. 
%Furthermore, we can easily construct an optimal 
%$1$-embedded graph on the torus that is not $1$-extendable as follows: 
%Consider such a graph $G$ such that $Q(G)$ has a quadrangulation $H$ as a subgraph, 
%and that every face of $H$ corresponds to an odd weighted region of $G$. 
%Then every non-crossing edge of $H$ is not extendable; 
%observe that $|V(H)|=C_o(G-V(H))$. 

Moreover, every O1PG and every O1PPG is not $3$-extendable since 
the graph has a vertex of degree exactly $6$. 
However, the property does not hold for optimal $1$-embedded graphs 
on the torus (or the Klein bottle); there exists infinitely many $8$-regular graphs whose quadrangular 
subgraphs are $4$-regular. 
At the end of the paper, we establish the following conjectures for those graphs: 

\begin{conj}\label{conj:1}
Let $G$ be an optimal $1$-embedded graph on the torus or the Klein bottle. 
If $Q(G)$ has no quadrangulation as a subgraph each of whose face corresponds 
to an odd weighted region of $G$, then $G$ is $1$-extendable.  
\end{conj}

\begin{conj}\label{conj2}
Every $8$-regular optimal $1$-embedded graph on the torus or the Klein bottle is 
$3$-extendable. 
\end{conj}

Note that for $2$-extendability, the statement for optimal $1$-embedded graphs on the 
torus and that for optimal $1$-embedded graphs on the Klein bottle might be substantially different; 
since the Klein bottle admits a separating simple closed curve that is not trivial, 
which is known as an {\em equator\/} of the surface.  

\bigskip

\noindent
{\bf Statements and Declarations}
\bigskip

This work was supported by JSPS KAKENHI Grant Number 23K03196. 
The authors have no relevant financial or non-financial interests to disclose.

%

%\bibliography{references} %.bibから拡張子を外した名前
%\bibliographystyle{elsarticle-num} %参考文献出力スタイル

% -----------------------------------------------------------------

\end{document}